\documentclass[11pt]{amsart}
\usepackage{upref,amssymb,eucal,epsf}

\usepackage{diagrams}
\diagramstyle[small,labelstyle=\scriptstyle,midshaft]
\newarrow{Same}=====
\newarrow{Corresponds}<--->
\newarrow{Dashonto}....{>>}
\newarrow{Dashembed}>...>

\numberwithin{equation}{section}

\theoremstyle{plain}
\newtheorem{theorem}{Theorem}[section]
\newtheorem{lemma}[theorem]{Lemma}
\newtheorem{proposition}[theorem]{Proposition}
\newtheorem{corollary}[theorem]{Corollary}

\newtheorem*{question}{Question}

\theoremstyle{definition}

\theoremstyle{remark}
\newtheorem{remark}[theorem]{Remark}
\newtheorem*{ack}{Acknowledgments}

\DeclareMathOperator{\Ext}{Ext}
\DeclareMathOperator{\shExt}{\mathcal E\mathit{xt}}

\DeclareMathOperator{\coker}{coker}
\DeclareMathOperator{\codim}{codim}

\DeclareMathOperator{\tr}{tr}

\DeclareMathOperator{\grade}{grade}

\DeclareMathOperator{\Gr}{Gr}

\DeclareMathOperator{\SL}{SL}
\DeclareMathOperator{\SO}{SO}
\DeclareMathOperator{\GL}{GL}

\DeclareMathOperator{\Span}{Span}

\DeclareMathOperator{\Pf}{Pf}
\DeclareMathOperator{\og}{\mathbb{OG}}
\DeclareMathOperator{\Sing}{Sing}
\DeclareMathOperator{\Char}{char}
\DeclareMathOperator{\coeff}{coeff}

\newcommand{\OG}[1]{\og_{2#1}}
\newcommand{\rest}[1]{{\mid}_{#1}}
\newcommand{\sm}[1]{\left(\begin{smallmatrix}#1\end{smallmatrix}\right)}

\newcommand{\mmu}{\boldsymbol{\mu}_2}
\newcommand{\aalpha}{\boldsymbol{\alpha}_2}

\begin{document}

\title[Non-Pfaffian Subcanonical Subschemes of Codimension $3$]
{Enriques Surfaces and other Non-Pfaffian 
Subcanonical Subschemes of Codimension $3$}

\author{David Eisenbud}
\address{Department of Mathematics\\ University of California,
Berkeley\\ Berkeley CA 94720 \\USA}
\email{de@msri.org}
\urladdr{http://www.msri.org/people/staff/de/}
\author{Sorin Popescu}
\address{Department of Mathematics\\ Columbia  University\\
New York, NY 10027\\ USA}
\email{psorin@math.columbia.edu} 
\urladdr{http://www.math.columbia.edu/\~{}psorin/}
\author{Charles Walter}
\address{Laboratoire J.-A.\ Dieudonn\'e (UMR 6621 du CNRS)\\
Universit\'e de Nice -- Sophia Antipolis\\ 06108 Nice Cedex 02\\
France}
\email{walter@math.unice.fr}
\urladdr{http://math1.unice.fr/\~{}walter/}

\thanks{Partial support for the authors during the preparation of
this work was provided by the NSF. The authors are also
grateful to MSRI Berkeley and the University of Nice Sophia-Antipolis
for their hospitality.}


\begin{abstract}
We give examples of subcanonical subvarieties of codimension $3$ in
projective $n$-space which are not {\em Pfaffian}, i.e.\ defined by
the ideal sheaf of submaximal Pfaffians of an alternating map of
vector bundles.  This gives a negative answer to a question asked by
Okonek \cite{Okonek}.

Walter \cite{Walter} had previously shown that a very large majority
of subcanonical subschemes of codimension $3$ in $\mathbb P^n$ are
Pfaffian, but he left open the question whether the exceptional
non-Pfaffian cases actually occur.  We give non-Pfaffian examples of
the principal types allowed by his theorem, including (Enriques)
surfaces in $\mathbb P^5$ in characteristic 2 and a smooth 4-fold in
$\mathbb P^7_{\mathbb C}$.

These examples are based on our previous work \cite{EPW} showing that
any strongly subcanonical subscheme of codimension $3$ of a Noetherian scheme
can be realized as a locus of degenerate intersection of a pair of
Lagrangian (maximal isotropic) subbundles of a twisted orthogonal
bundle.
\end{abstract}

\maketitle

\label{introduction}

There are many relations between the vector bundles $\mathcal E$ on a
nonsingular algebraic variety $X$ and the subvarieties $Z \subset X$.
For example, a globalized form of the Hilbert-Burch theorem allows one
to realize any codimension $2$ locally Cohen-Macaulay subvariety as a
degeneracy locus of a map of vector bundles.  In addition the Serre
construction gives a realization of any subcanonical codimension $2$
subvariety $Z \subset X$ as the zero locus of a section of a rank $2$
vector bundle on $X$.

The situation in codimension $3$ is more complicated.  In the local
setting, Buchsbaum and Eisenbud \cite{BE} described the structure of
the minimal free resolution of a Gorenstein (i.e.\ subcanonical)
codimension $3$ quotient ring of a regular local ring.  Their
construction can be globalized: If $\phi: \mathcal E \to \mathcal
E^*(L)$ is an alternating map from a vector bundle $\mathcal E$ of
\textbf{odd} rank $2n+1$ to a twist of its dual by a line bundle $L$,
then the $2n\times 2n$ Pfaffians of $\phi$ define a degeneracy locus
$Z = \{ x \in X \mid \dim \ker(\phi(x)) \geq 3 \}$.  If $X$ is
nonsingular and $\codim(Z)=3$, the largest possible value, (or more
generally if $X$ is locally Noetherian and $\grade(Z)=3$), then, after
twisting $\mathcal E$ and $L$ appropriately, $\mathcal O_Z$ has a
symmetric locally free resolution
\begin{equation}
\label{Okonek.resol}
0 \to L \xrightarrow{p^*} \mathcal E \xrightarrow{\phi} \mathcal
E^*(L) \xrightarrow{\,p\,} \mathcal O_X \to \mathcal O_Z \to 0
\end{equation}
with $p$ locally the vector of submaximal Pfaffians of $\phi$ and with
$Z$ subcanonical with $\omega_Z \cong \omega_X (L^{-1}) \rest{Z}$.
(The reader will find a general discussion of {\em subcanonical}
subschemes in the introduction of our paper \cite{EPW} and in Section
\ref{review} below; for the purpose of this introduction it may
suffice to know that a codimension 3 subvariety $Z$ of $\mathbb P^n$
is {\em subcanonical} if it is (locally) Gorenstein with canonical
line bundle $\omega_Z =\mathcal O_Z(d)$.)

Okonek \cite{Okonek}
called such a $Z$ a {\em Pfaffian subvariety}, and he asked:

\begin{question}[Okonek]
Is every smooth subcanonical subvariety of codimension $3$ in a smooth
projective variety a Pfaffian subvariety\textup?
\end{question}

This paper is one of a series which we have written in response to
Okonek's question.  In the first paper (Walter \cite{Walter}) one of
us found necessary and sufficient numerical conditions for a
codimension $3$ subcanonical subscheme $X \subset \mathbb P^N_k$ to be
Pfaffian, but was unable to give subcanonical schemes failing the
numerical conditions.  In the second paper \cite{EPW} we found a
construction for codimension $3$ subcanonical subschemes which is more
general than that studied by Okonek, and proved that this construction
gives all subcanonical subschemes of codimension 3 satisfying a
certain necessary lifting condition.  (This lifting condition always
holds if the ambient space is $\mathbb P^n$ or a Grassmannian of
dimension at least $4$.)

The construction is as follows: Let $\mathcal E, \mathcal F$ be a pair
of Lagrangian subbundles of a twisted orthogonal bundle with
$\dim_{k(x)} \bigl [ \mathcal E(x) \cap \mathcal F(x)\bigr]$ always
odd. There is a natural scheme structure on the degeneracy
locus
\begin{equation}
\label{Z}
Z := \{ x \in X \mid \dim_{k(x)} \bigl[ \mathcal E(x) \cap \mathcal F(x)
\bigr] \geq 3 \}
\end{equation}
which is subcanonical if $\grade(Z) = 3$; the structure theorem
asserts that every strongly subcanonical codimension 3 subscheme (in
the sense of the definition in the introduction of our paper
\cite{EPW}) arises in this way.

Another more complicated but also more explicit description of the
structure theorem may be given as follows (see our paper \cite{EPW}
for more details): If one has a subcanonical subvariety with $\omega_Z
\cong L \rest Z$, then the isomorphism $\mathcal O_Z \cong
\omega_Z(L^{-1}) \rest Z \cong \shExt^3_{\mathcal O_X}(\mathcal O_Z,
\omega_X(L^{-1}))$ is a sort of symmetry in the derived category.  Let
$L_1 := \omega_X(L^{-1})$.  Then $\mathcal O_Z$ should have a locally
free resolution which is symmetric in the derived category.  This is
the case if $\mathcal O_Z$ has a locally free resolution with a
symmetric quasi-isomorphism into the twisted shifted dual complex:
\begin{diagram}[LaTeXeqno]
\label{intro.diag}
0 & \rTo & L_1 & \rTo & \mathcal E & \rTo ^\psi & \mathcal G & \rTo &
\mathcal O_X & \rTo & \mathcal O_Z & \rTo & 0 \\
&& \dSame && \dTo <\phi && \dTo >{\phi^*} && \dSame && \dTo >\eta
<\cong \\
0 & \rTo & L_1 & \rTo & \mathcal G^*(L_1) & \rTo ^{-\psi^*} & \mathcal
E^*(L_1)& \rTo & \mathcal O_X & \rTo & \shExt^3_{\mathcal O_X} (\mathcal
O_Z, L_1) & \rTo & 0
\end{diagram}
%
If $\phi^*$ is the identity, then we get the Pfaffian resolution of
\eqref{Okonek.resol}.  But in general it is not.  However, the
resolutions in \eqref{intro.diag} mean that $Z$ can be interpreted as
the degeneracy locus of \eqref{Z} for the Lagrangian subbundles
$\mathcal E,\mathcal G^*(L) \subset \mathcal G \oplus \mathcal
G^*(L)$.

In this paper we give examples of non-Pfaffian codimension $3$
subcanonical subvarieties.  Thus we give an answer ``No'' to Okonek's
question.  Because of the numerical conditions of \cite{Walter}, our
examples fall into four types:
\begin{enumerate}
\item Surfaces $S \subset \mathbb P^5_k$ with $\Char(k)=2$ and with
$\omega_S \cong \mathcal O_S(2d)$ such that $h^1(\mathcal O_S(d))$ is
odd.
\item Fourfolds $Y \subset \mathbb P^7$ with $\omega_Y \cong \mathcal
O_Y(2d)$ such that $h^2(\mathcal O_Y(d))$ is odd.
\item Fourfolds $Z \subset \mathbb P^7_{\mathbb R}$ such that
$\omega_Z \cong \mathcal O_Z(2d)$ such that the cup product pairing on
$H^2(\mathcal O_Z(d))$ is positive definite (or otherwise not
hyperbolic).
\item Examples in ambient varieties other than $\mathbb P^N$.
\end{enumerate}
Our examples include among others: codimension $3$ Schubert
subvarieties of orthogonal Grassmannians, non-Pfaffian fourfolds of
degree $336$ in $\mathbb P^7$, and nonclassical Enriques surfaces in
$\mathbb P^5$ in characteristic $2$.

\subsection{Structure of the paper}

In \S\ref{sect.grass} we review well known about quadratic forms,
orthogonal Grassmannians and Pfaffian line bundles.  In \S
\ref{review} and \S \ref{sect.Pfaffian} we review the machinery and
results of our previous papers \cite{EPW} and \cite{Walter}.  The rest
of the paper is devoted to the examples of non-Pfaffian subcanonical
subschemes.  The first example, in \S \ref{sect.Schubert}, is a
codimension $3$ Schubert variety in $\OG{n}$.  In \S \ref{sect.P7}
we construct a smooth non-Pfaffian subcanonical fourfold of degree
$336$ in $\mathbb P^7$ with $K = 12H$.  In \S \ref{sect.real} we
construct some subcanonical fourfolds in $\mathbb P^7_{\mathbb R}$
which are not Pfaffian over $\mathbb R$ but become Pfaffian over
$\mathbb C$.  In \S \ref{sect.rp2} we analyze Reisner's example in
characteristic $2$ of the union of $10$ coordinate $2$-planes in
$\mathbb P^5$.

In the last two sections \S\S \ref{sect.Enriques},\ref{sect.Frob} we
analyze nonclassical Enriques surfaces in characteristic $2$ and their
Fano models.  As a consequence we obtain a (partial) moduli
description of Fano-polarized unnodal nonclassical Enriques surfaces
as a ``quotient'' of $\og_{20}$ modulo $\SL_6$ (Corollary
\ref{moduli}).  In \S \ref{sect.Frob} we identify the closure of the
locus of $\aalpha$-surfaces within $\og_{20}$ by calculating the
action of Frobenius on our resolutions.

The philosophy that (skew)-symmetric sheaves should have locally free
resolutions that are (skew)-symmetric up to quasi-isomorphism is also
pursued in \cite{EPW3} and \cite{Walter2}.  The former deals primarily
with methods for constructing explicit locally free resolutions for
(skew)-symmetric sheaves on $\mathbb P^n$.  The latter investigates
the obstructions (in Balmer's derived Witt groups \cite{Balmer}) for
the existence of a genuinely (skew)-symmetric resolution.

\begin{ack}
We are grateful to Igor Dolgachev, Andr\'e Hirschowitz, and N.\
Shepherd-Barron for many useful discussions or suggestions.

The second and third authors would like to thank the Mathematical
Sciences Research Institute in Berkeley for its support while part of
this paper was being written.  The first and second authors also thank
the University of Nice--Sophia Antipolis for its hospitality.  

Many of the diagrams were set using the {\tt diagrams.tex}
package of Paul Taylor.
\end{ack}

\section{Orthogonal Grassmannians}
\label{sect.grass}

In this section we bring together a number of facts about orthogonal
Grassmannians which we will need throughout the paper.  All the
results are well known, so we omit proofs.  

\subsection{Quadratic forms \cite{Knus} \cite{Pfister}}
Let $V$ be vector space of even dimension $2n$ over a field $k$.  (We
impose no restrictions on $k$; it need not be algebraically closed nor
of characteristic $\neq 2$.)  A quadratic form on $V$ is a homogeneous
quadratic polynomial in the linear forms on $V$, i.e.\ a member $q \in
S^2(V^*)$.  The symmetric bilinear form $b: V\times V \to k$
associated to $q$ is given by the formula
\begin{equation}
\label{ass.bilin.form}
b(x,y) := q(x+y) - q(x) - q(y).
\end{equation}
The quadratic form $q$ is {\em nondegenerate} if $b$ is a perfect
pairing.

A {\em Lagrangian subspace} $L \subset (V,q)$ is a subspace of
dimension $n$ such that $q \rest{L} \equiv 0$.  There exist such
subspaces if and only if the quadratic form $q$ is {\em hyperbolic}.
By definition, this means that there exists a coordinate system on $V$
in which one can write $q = \sum_{i=1}^n x_i x_{n+i}$.  Lagrangian
subspaces $L$ satisfy $L= L^\perp$, and the converse is true if
$\Char(k) \neq 2$.

If $X$ is an algebraic variety over the field $k$, then the quadratic
form $q$ makes $V\otimes \mathcal O_X$ into an orthogonal bundle.  A
{\em Lagrangian subbundle} $\mathcal E \subset V\otimes \mathcal O_X$
is a subbundle of rank $n$ such that all fibers $\mathcal E(x) \subset
V$ are Lagrangian with respect to $q$.

\subsection{Orthogonal Grassmannians \cite{Mukai} \cite{Mumford}
\cite{Pressley}}
Let $V$ be a vector space of even dimension $2n$ with the hyperbolic
quadratic form $q=\sum_{i=1}^n x_i x_{n+i}$.

The Lagrangian subspaces of $(V,q)$ form two disjoint families such
that $\dim(L \cap W) \equiv n \pmod 2$ if $L$ and $W$ are Lagrangian
subspaces in the same family, while $\dim(L \cap W) \equiv n+1 \pmod
2$ otherwise.

Each family of Lagrangian subspaces is parametrized by the
$k$-rational points of an orthogonal Grassmannian (or spinor variety)
$\OG{n}$ of dimension $n(n-1)/2$.  Let $U$ (resp.\ $W$) be a
Lagrangian subspace such that $\dim(L\cap U)$ is even (resp.\ $\dim(L
\cap W)$ is odd) for all $L$ in the first family.  Then $\OG{n}$
contains the following Schubert varieties:
\begin{align}
\sigma_{2}(U) & := \{ L \mid \dim(L\cap U) \geq 2 \}, \\
\sigma_{3}(W) & := \{ L \mid \dim(L\cap W) \geq 3 \}, \\
\sigma_{5}(W) & := \{ L \mid \dim(L\cap W) \geq 5 \}. 
\end{align}
We have the following basic results about these subvarieties:

\begin{lemma}
\label{cells}
\textup{(a)} The subvariety $\sigma_2(U) \subset \OG{n}$ is of
codimension $1$, and its complement \textup(parametrizing Lagrangian
subspaces complementary to $U$\textup) is an affine space $\mathbb
A^{n(n-1)/2}_k$.  

\textup{(b)} The subvarieties $\sigma_3(W)$ and $\sigma_5(W)$ are of
codimensions $3$ and $10$ in $\OG{n}$, respectively.  Moreover,
$\sigma_3(W)$ is nonsingular outside $\sigma_5(W)$.
\end{lemma}

The reason behind Lemma \ref{cells}(a) is that if one fixes a
Lagrangian subspace $U^*$ complementary to $U$, then the Lagrangian
subspaces complementary to $U$ are the graphs of alternating maps $U
\to U^*$.  This fact also underlies the following lemma.

\begin{lemma}
\label{tangent}
Suppose $L \subset (V,q)$ is a Lagrangian subspace.  Then there is a
natural identification $T_{[L]} (\OG{n}) \cong \Lambda^2 L^*$.
Moreover, if $\dim(L \cap W) = 3$, then there is a natural
identification
\[
T_{[L]} \left( \sigma_3(W) \right) \cong \{ f \in \Lambda^2 L^* \mid f
\rest{L \cap W} \equiv 0 \}.
\]
\end{lemma}

\subsection{Pfaffian line bundles \cite{Mukai} \cite{Pressley}}

The Picard group of $\OG{n}$ is isomorphic to $\mathbb Z$.  Its
positive generator $\mathcal O(1)$ is the {\em universal Pfaffian line
bundle}.  The divisors $\sigma_2(U)$ of Lemma \ref{cells}(a)
are zero loci of particular sections of $\mathcal O(1)$.  The
canonical line bundle of $\OG{n}$ is $\mathcal O(2-2n)$.

On the orthogonal Grassmannian $\OG{n}$ there is a universal
short exact sequence of vector bundles
\begin{equation}
\label{univ.ex.seq}
0\to \mathcal S\to V\otimes \mathcal O_{\OG{n}}\to \mathcal S^*\to 0
\end{equation}
with $\mathcal S$ the {\em universal Lagrangian subbundle} of rank
$n$.  The determinant line bundle related to the Pl\"ucker embedding
is $\det(\mathcal S^*) = \mathcal O(2)$.

The universal Lagrangian subbundle has the following universal
property: If $X$ is an algebraic variety and $\mathcal E\subset
V\otimes \mathcal O_X$ is a Lagrangian subbundle over $X$ whose fibers
lie in the family parametrized by $\OG{n}$, then there exists a unique
map $f: X\to \OG{n}$ such that the natural exact sequence $0 \to
\mathcal E \to V \otimes \mathcal O_X \to \mathcal E^* \to 0$ is the
pullback along $f$ of the universal exact sequence
\eqref{univ.ex.seq}.  The {\em Pfaffian line bundle} of $\mathcal E$
is then $\Pf (\mathcal E) := f^* \bigl( \mathcal O(-1) \bigr)$.  It is
a canonically defined line bundle on $X$ such that $\Pf(\mathcal E)
^{\otimes 2} \cong \det(\mathcal E)$.

\section{The Lagrangian subbundle construction}
\label{review}

Let $X$ be a nonsingular algebraic variety over a field $k$, and let
$Z\subset X$ be a codimension $3$ subvariety.  In \cite{EPW} we called
$Z\subset X$ {\em subcanonical} if the sheaf $\mathcal O_Z$ has local
projective dimension over $\mathcal O_X$ equal to its grade (in this
case $3$), and the relative canonical sheaf $\omega_{Z/X}:=
\shExt^3_{\mathcal O_X}(\mathcal O_Z,\mathcal O_X)$ is isomorphic to
the restriction to $Z$ of a line bundle $L$ on $X$.  We called $Z$
{\em strongly subcanonical} if in addition the map
\[
\Ext^3_{\mathcal O_X}(\mathcal O_Z,L^{-1})\rightarrow
\Ext^3_{\mathcal O_X}(\mathcal O_X,L^{-1})=H^3(X,L^{-1}).
\]
induced by the surjection $\mathcal O_X\to \mathcal O_Z$ is zero.
The last condition is immediate if $X=\mathbb P^n$, for $n\ge 4$,
as the right cohomology group is always zero.  Thus all subcanonical
subvarieties of $\mathbb P^n$ are strongly subcanonical.

The method we use to construct subcanonical subvarieties of codimension
$3$ was developed in \cite{EPW}.  In fact, all our examples can be
constructed using the following result, which is a special case of
\cite{EPW} Theorem 3.1.

\begin{theorem}
\label{split}
Let $V$ be a vector space of even dimension $2n$ over a field $k$, and
$q = \sum_{i=1}^n x_i x_{n+i}$ a hyperbolic quadratic form on $V$.
Suppose that $X$ is a nonsingular algebraic variety over $k$, and that
$\mathcal E \subset V \otimes \mathcal O_X$ is a Lagrangian subbundle
with Pfaffian line bundle $L := \Pf(\mathcal E)$ satisfying
$L^{\otimes 2} \cong \det(\mathcal E)$.  Let $W \subset (V,q)$ be a
Lagrangian subspace such that $\dim_{k(x)} \left[ \mathcal E(x) \cap W
\right]$ is odd for all $x$.  If $W$ is sufficiently general, then
\begin{equation}
Z_W = \{ x \in X \mid \dim_{k(x)} \left( \mathcal E(x) \cap W \right)
\geq 3 \}, 
\end{equation}
is a \textup(strongly\textup) subcanonical subvariety of codimension
$3$ in $X$ with
\[
\omega_{Z_W} \cong \left( \omega_X \otimes \det(\mathcal E^*) \right)
\rest {Z_W}
\]
and with symmetrically quasi-isomorphic locally free resolutions
\begin{footnotesize}
\begin{diagram}[LaTeXeqno]
\label{dual.diag}
0 & \rTo & L^{\otimes 2} & \rTo & \mathcal E(L) & \rTo ^\psi & W^*
\otimes \mathcal O_X(L) & \rTo & \mathcal O_X & \rOnto & \mathcal O_{Z_W}
\\
&& \dSame && \dTo <\phi && \dTo >{\phi^*} && \dSame && \dTo >\eta
<\cong \\
0 & \rTo & L^{\otimes 2} & \rTo & W \otimes \mathcal O_X(L) & \rTo
^{-\psi^*} & \mathcal E^*(L)& \rTo & \mathcal O_X & \rOnto &
\shExt^3_{\mathcal O_X} (\mathcal O_{Z_W}, L^{\otimes 2})
\end{diagram}%
\end{footnotesize}%
Moreover, if $W$ is sufficiently general and $\Char(k) = 0$, then
$\Sing({Z_W})$ is of codimension $10$ in $X$.
\end{theorem}

According to Fulton-Pragacz \cite{FP} (6.5), the fundamental class of
${Z_W}$ in the Chow ring of $X$ is $\frac 14 (c_1 c_2 - 2 c_3)$, where
$c_i := c_i(\mathcal E^*)$.

To check the smoothness of ${Z_W}$ in characteristic $0$ we proceed as
follows.  The universal property of the orthogonal Grassmannian gives
us a morphism $f: X \to \OG{n}$.  The group $\SO_{2n}$ acts on
$\OG{n}$, and it translates the Schubert varieties into other Schubert
varieties: $g \cdot \sigma_i(\Lambda) = \sigma_i(g\Lambda)$.  Since
the characteristic is $0$, Kleiman's theorem on the transversality of
the general translate \cite{Kleiman} applies.  So there exists a
nonempty Zariski open subset $U \subset \SO_{2n}$ such that if $g \in
U$, then $\sigma_5(g\Lambda)$ intersects $f(X)$ transversally in
codimension $10$, and $\sigma_3(g\Lambda)$ intersects $f(X)$
transversally in codimension $3$.  Moreover, since $\SO_{2n}$ is a
rational variety and $k$ is an infinite field, $U$ contains
$k$-rational points.  So if $W := g\Lambda$ with $g$ a $k$-rational
point in $U$, then $W \subset (V,q)$ is a Lagrangian subspace defined
over $k$ such that ${Z_W} := f^{-1}(\sigma_3(W))$ is of codimension
$3$ and is smooth outside $f^{-1}(\sigma_5(W))$, which is of
codimension $10$ in $X$.

In characteristic $p$, it will be more complicated to prove that
${Z_W}$ is smooth (see Lemma \ref{transversal} below).

\section{Pfaffian subschemes}
\label{sect.Pfaffian}

Theorem \ref{split} allows us to construct subcanonical subvarieties
$Z \subset X$ of codimension $3$ with resolutions which are not
symmetric.  Such a $Z$ might nevertheless possess other locally free
resolutions which are symmetric (i.e.\ Pfaffian).  However, the third
author gave in \cite{Walter} a necessary and sufficient condition for
a subcanonical subscheme of codimension $3$ in projective space to be
Pfaffian.

So suppose $Z \subset \mathbb P^{n+3}$ is a subcanonical subscheme of
codimension $3$ and dimension $n$ such that $\omega_Z \cong \mathcal
O_Z(\ell)$.  If $n$ and $\ell$ are both even, then Serre duality
defines a nondegenerate bilinear form
\begin{equation}
\label{cup}
H^{n/2}(\mathcal O_Z(\ell/2)) \times H^{n/2}(\mathcal O_Z(\ell/2))
\xrightarrow {\cup} H^n(\mathcal O_Z(\ell)) \cong k.
\end{equation}
This bilinear form is symmetric if $n\equiv 0 \pmod 4$, and it is
skew-symmetric if $n \equiv 2 \pmod 4$.  In characteristic $2$ it is
symmetric in both cases.  The result proven was the following.

\begin{theorem}[Walter \cite{Walter}]
\label{criterion}
Suppose that $Z \subset \mathbb P^{n+3}$ is a subcanonical subscheme
of codimension $3$ and dimension $n \geq 1$ over the field $k$ such
that $\omega_Z \cong \mathcal O_Z (\ell)$.  Then $Z$ is Pfaffian if
and only at least one of the following conditions holds\textup:
\textup{(i)} $n$ or $\ell$ is odd, \textup{(ii)} $n \equiv 2 \pmod 4$
and $\Char(k) \neq 2$, or \textup{(iii)} $n$ and $\ell$ are even, and
$H^{n/2}(\mathcal O_Z(\ell/2))$ is of even dimension and contains a
subspace which is Lagrangian with respect to the cup product pairing
of \eqref{cup}.
\end{theorem}

Another way of stating the theorem is that a subcanonical subscheme $Z
\subset \mathbb P^{n+3}$ of dimension $n$ with $\omega_Z \cong
\mathcal O_Z(\ell)$ is not Pfaffian if and only either
\begin{enumerate}
\item [(a)]
$n \equiv 0 \pmod 4$, or else
\item [(b)] 
$n\equiv 2 \pmod 4$ and $\Char(k) = 2$;
\end{enumerate}%
and at the same time either
\begin{enumerate}
\item [(c)] $\ell$ is even and $H^{n/2}(\mathcal O_X(\ell/2))$ is
  odd-dimensional, or else
\item [(d)]
$\ell$ is even and $H^{n/2}(\mathcal O_X(\ell/2))$ is even-dimensional
but has no subspace which is Lagrangian with respect to the cup
product pairing of \eqref{cup}.
\end{enumerate}

The simplest cases where non-Pfaffian subcanonical subschemes of
dimension $n$ in $\mathbb P^{n+3}_k$ could exist are therefore

\begin{enumerate}
\item Surfaces $S \subset \mathbb P^5_k$ with $\Char(k)=2$ and with
$\omega_S \cong \mathcal O_S(2d)$ such that $h^1(\mathcal O_S(d))$ is
odd.
\item Fourfolds $Y \subset \mathbb P^7$ with $\omega_Y \cong \mathcal
O_Y(2d)$ such that $h^2(\mathcal O_Y(d))$ is odd.
\item Fourfolds $Z \subset \mathbb P^7_{\mathbb R}$ such that
$\omega_Z \cong \mathcal O_Z(2d)$ such that the cup product pairing on
$H^2(\mathcal O_Z(d))$ is positive definite (or otherwise not
hyperbolic).
\end{enumerate}

\section{The codimension $3$ Schubert variety in $\OG{n}$}
\label{sect.Schubert}

Our first examples of non-Pfaffian subcanonical subschemes of
codimension $3$ are the codimension $3$ Schubert varieties of
$\OG{n}$.

\begin{theorem}
\label{Schubert}
Let $n \geq 4$, and let $\sigma_3(W) \subset \OG{n}$ be one of the
codimension $3$ Schubert varieties of Lemma \ref{cells}\textup{(b)}.
Then $Z$ is \textup(strongly\textup) subcanonical with $\omega_Z \cong
\mathcal O_Z (4-2n)$ but is not Pfaffian.
\end{theorem}

\begin{proof}
We apply Theorem \ref{split} with $X := \OG{n}$ and $\mathcal E :=
\mathcal S$ the universal Lagrangian subbundle.  The degeneracy locus
is $Z := \sigma_3(W)$ and it has resolutions
\begin{footnotesize}
\begin{diagram}[LaTeXeqno]
\label{Schubert.diag}
0 & \rTo & \mathcal O(-2) & \rTo & \mathcal S(-1) & \rTo ^\psi & W^*
\otimes \mathcal O(-1) & \rTo & \mathcal O & \rOnto & \mathcal O_Z
\\
&& \dSame && \dTo <\phi && \dTo >{\phi^*} && \dSame && \dTo >\eta
<\cong \\
0 & \rTo & \mathcal O(-2) & \rTo & W \otimes \mathcal O(-1) & \rTo
^{-\psi^*} & \mathcal S^*(-1)& \rTo & \mathcal O & \rOnto &
\shExt^3 (\mathcal O_Z, \mathcal O(-2))
\end{diagram}%
\end{footnotesize}%
Moreover, $\omega_Z \cong \omega_{\OG{n}}(2) \rest{Z} \cong \mathcal
O_Z(4-2n)$ since $\omega_{\OG{n}} \cong \mathcal O(2-2n)$.

It remains to show that $Z$ is not Pfaffian.  Let $U$ be a general
totally isotropic subspace of dimension $n-3$, and let $i :
\og_{6}(U^\perp/U) \hookrightarrow \OG{n}$ be the natural embedding.
If $Z$ were Pfaffian, its symmetric resolution would be of the form
\begin{equation}
\label{Schubert.Pfaff}
0 \to \mathcal O(-2) \to \mathcal E(-1) \to \mathcal E^*(-1) \to
\mathcal O \to \mathcal O_Z \to 0
\end{equation}
(cf.\ \eqref{Okonek.resol} in the introduction).  It would thus pull
back to a symmetric resolution of $i^{-1}(Z)$ on $\og_{6}(U^\perp/U)$
of the same form.

Now $\og_{6}$ parametrizes one family of $\mathbb P^2$'s contained in
a smooth hyperquadric in $\mathbb P^5$.  Thus $\og_{6} \cong \mathbb
P^3$ (see Jessop \cite{Jessop} or Griffiths-Harris \cite{GH}), and
$i^{-1}(Z)$ is a codimension $3$ Schubert subvariety of $\og_{6} \cong
\mathbb P^3$ and indeed a point $Q$.  If we pull the resolution
\eqref{Schubert.Pfaff} back to $\mathbb P^3$ and twist, we get a
resolution
\[
0 \to \mathcal O_{\mathbb P^3} (-3) \to i^* \mathcal E(-2) \to i^*
\mathcal E^*(-2) \to \mathcal O_{\mathbb P^3}(-1) \to \mathcal O_Q(-1)
\to 0.
\]
This gives $1 = \chi(\mathcal O_Q(-1)) = 2 \chi(i^* \mathcal E(-2))$,
which is a contradiction.  This completes the proof.
\end{proof}

The contradiction obtained at the end of the proof above is of the
same nature as the obstructions in Theorem \ref{criterion} (see
\cite{Walter}).  The contradiction can also be derived from \cite{EPW}
Theorem 7.2.

\section{A non-Pfaffian subcanonical $4$-fold in $\mathbb P^7$}
\label{sect.P7}

In this section we use the construction of Theorem \ref{split} to give
an example of a smooth non-Pfaffian subcanonical $4$-fold in $\mathbb
P^7$ over an arbitrary infinite field $k$.

Let $V = H^0(\mathcal O_{\mathbb P^7}(1))^*$, so that $\mathbb P^7$ is
the projective space of lines in $V$.  We fix an identification
$\Lambda^8 V \cong k$.  The $70$-dimensional vector space $\Lambda^4
V$ has the quadratic form $q(u) = u^{(2)}$, the divided square in the
exterior algebra.  (If $\Char(k) \neq 2$, then $u^{(2)} = \frac 12 u
\wedge u$.)  The associated symmetric bilinear form is $b(u,v) = u
\wedge v$.  This quadratic form is hyperbolic because it has
Lagrangian subspaces, as we shall see in a moment.

We now look for a Lagrangian subbundle of $\Lambda^4 V \otimes
\mathcal O_{\mathbb P^7}$.  On $\mathbb P^7$ there is a canonical map
$\mathcal O_{\mathbb P^7} \to V \otimes \mathcal O_{\mathbb P^7}(1)$.
Using this, one constructs a Koszul complex of which one part is:
\begin{equation}
\label{Koszul}
\dotsb \to \Lambda^3 V \otimes \mathcal O_{\mathbb P^7}(-1)
\xrightarrow{\,d\,} \Lambda^4 V \otimes \mathcal O_{\mathbb P^7} \to
\Lambda^5 V \otimes \mathcal O_{\mathbb P^7}(1) \to \dotsb
\end{equation}
The image and cokernel of $d$ fit into an exact sequence
\begin{equation}
\label{Euler}
0 \to \Omega^4_{\mathbb P^7}(4) \to \Lambda^4 V \otimes \mathcal
O_{\mathbb P^7} \to \Omega^3_{\mathbb P^7}(4) \to 0.
\end{equation}
The fiber of $\Omega^4_{\mathbb P^7}(4)$ over a point $\overline \xi
\in \mathbb P^7$, viewed as a subspace of $\Lambda^4 V$, is the image
over $\overline \xi$ of the map $d$ of \eqref{Koszul}, namely $W_\xi
:= \xi \wedge \Lambda^3 V$.  These spaces are all totally isotropic of
dimension $35$, half the dimension of $\Lambda^4 V$.  Hence
$\Omega^4_{\mathbb P^7}(4)$ is a Lagrangian subbundle of $\Lambda^4 V
\otimes \mathcal O_{\mathbb P^7}$.

The fibers $W_\xi$ of $\Omega^4_{\mathbb P^7}(4)$ are all members of
the same family of Lagrangian subspaces of the hyperbolic quadratic
space $\Lambda^4 V$.  Let $W$ be a general member of that family, and
(to abuse notation) $W^*$ a general member of the opposite family of
Lagrangian subspaces of $\Lambda^4 V$.  Then $\Lambda^4 V = W \oplus
W^*$, and the symmetric bilinear form on $\Lambda^4 V$ induces a
natural isomorphism between $W^*$ and the dual of $W$.  Finally, let
\[
Z_W := \{ \overline \xi \in \mathbb P^7 \mid \dim (W_\xi \cap W) \geq
3 \}.
\]
We now apply Theorem \ref{split} and obtain the following result.

\begin{theorem}
\label{counterexample}
If the base field $k$ is infinite, and if $W \subset \Lambda^4 V$ is a
Lagrangian subspace in the same family as the $W_\xi$ which is
sufficiently general, then $Z_W$ is a nonsingular subvariety of
$\mathbb P^7$ of dimension $4$ which is subcanonical but not Pfaffian.
It is of degree $336$ and has $\omega_{Z_W} \cong \mathcal
O_{Z_W}(12)$, and its structural sheaf has locally free resolutions
\textup(where we write $\mathcal O := \mathcal O_{\mathbb
P^7}$\textup)
\begin{small}
\begin{diagram}[LaTeXeqno]
\label{Omegas}
0 & \rTo & \mathcal O(-20) & \rTo & \Omega^4_{\mathbb
P^7}(-6) & \rTo^\psi & W^* \otimes \mathcal O(-10) &
\rTo & \mathcal O & \rOnto & \mathcal O_{Z_W} \\
&& \dSame && \dTo <\phi && \dTo >{\phi^*} && \dSame && \dTo >\eta
<\cong \\
0 & \rTo & \mathcal O(-20) & \rTo & W \otimes \mathcal
O(-10) & \rTo^{-\psi^*} & \Omega^3_{\mathbb P^7}(-6) &
\rTo & \mathcal O & \rOnto & \omega_{Z_W}(-12)
\end{diagram}
\end{small}
\end{theorem}

\begin{proof}
The existence of $Z_W$ and the resolutions follow from Theorem
\ref{split}, using the fact that $\det(\Omega^4_{\mathbb P^7}(4)) =
\mathcal O_{\mathbb P^7}(-20)$.  The subcanonical subvariety $Z_W$ is
not Pfaffian because $\omega_Z = \mathcal O_Z(12)$ and $h^2(\mathcal
O_Z(6)) = 1$ (as one sees from \eqref{Omegas}).  Using the resolutions
and a computer algebra package, one computes the Hilbert polynomial of
$Z_W$ as:
\[
336 \binom{n+3}4 - 2520 \binom{n+2}3 + 9814 \binom{n+1}2 -25571n
+49549
\]
So the degree of $Z_W$ is $336$.  In characteristic $0$ the smoothness
of the general $Z_W$ follows from Kleiman's transversality theorem as
explained after Theorem \ref{split}.  In characteristic $p$ the
smoothness follows from Lemma \ref{transversal} below.
\end{proof}

The free resolution of $R/I_{Z_W}$ is of the form
\begin{multline}
\label{free.res.P7}
0 \to R(-14) \to R(-13)^8 \to R(-12)^{28} \oplus R(-20) \\ \to
R(-11)^{56} \to R(-10)^{35} \to R \to R/I_{Z_W} \to 0.
\end{multline}

\begin{lemma}
\label{transversal}
In Theorem \ref{counterexample}, if $k$ is an infinite field and $W$
is general, then $Z_W$ is smooth.
\end{lemma}

\begin{proof}
Let
\begin{multline*}
\Upsilon:= \{ W \in \og_{70} \mid \mathbb P^7 \cap \sigma_5(W) =
\varnothing, \\ \text{ and } \mathbb P^7 \cap \sigma_3(W) = Z_W
\text{ is nonsingular of dimension $4$} \}
\end{multline*}

We will show that $\Upsilon$ is nonempty without invoking Kleiman's
theorem on the transversality of a general translate.  We will do this
by showing that the complement has dimension at most $\dim
\og_{70}-1$.

Let $n:=35$, and let $\Gr(3,2n)$ be the Grassmannian parametrizing
three-dimensional subspaces of $\Lambda^4 V$.  Let 
\[
Y := \{ (\xi,F,W) \in \mathbb P^7 \times \Gr(3,2n) \times \OG{n}
\mid F \subset W_\xi \cap W \}.
\]
To choose a point of $Y$, one first chooses $\xi \in \mathbb P^7$,
then $F$ in a fiber of a Grassmannian bundle with fiber $\Gr(3,n)$ ,
then $W$ in a fiber of an orthogonal Grassmannian bundle with fiber
$\OG{(n-3)}$.  We see that $Y$ is nonsingular with $\dim(Y) = \dim
\OG{n} + 4$.  A similar computation shows that
\[
Y_1 := \{ (\xi,F,W) \in Y \mid \dim(W_\xi \cap W) \geq 5 \}
\]
is of dimension equal to $\dim \OG{n}-3$.  So the image of $Y_1$
under the projection $\pi: Y \to \OG{n}$ is not dominant.  Over a
point $W \notin \pi(Y_1)$, one has $\pi^{-1}(W) = Z_W$.  Hence it
is enough to show that the dimension of
\[
Y_2 := \{ (\xi,F,W) \in Y \mid  d\pi : T_{(\xi,F,W)}Y \to T_W \OG{n}
\text{ is not surjective} \}
\]
is at most $\dim \OG{n} -1$.  The essential question is: if
$(\xi,F,W) \in Y - Y_1$, so that $\pi^{-1}(W) = Z_W$, when is the
dimension of $T_\xi Z_W = \ker (d\pi)_{(\xi,F,W)}$ equal to $4$, and
when is it more?

Now $T_\xi Z_W = T_{W_\xi}(\mathbb P^7 \cap \sigma_3(W))$.  In order
to identify this space we must identify $T_{W_\xi}(\OG{n})$ and the
two subspaces $T_{W_\xi}\mathbb P^7$ and $T_{W_\xi}(\sigma_3(W))$.

There is a well known isomorphism $T_{W_\xi}(\OG{n}) \cong \Lambda^2
W_\xi^*$.  Essentially, given a first-order deformation of $W_\xi$ as
a Lagrangian subspace of $\Lambda^4 V$, any vector $w \in W_\xi$
deforms within the deforming subspace as a $w + \varepsilon u(w)$ with
$u(w)$ well-defined modulo $W_\xi$.  This gives a map $W_\xi \to
\Lambda^4 V/W_\xi \cong W_\xi^*$.  The vector $w+\varepsilon u(w)$
remains isotropic if and only if $u$ is alternating.  To describe $u$
as an alternating bilinear form on $W_\xi$, one looks at
\[
\langle w_1 , w_2 + \varepsilon u(w_2) \rangle = \varepsilon \langle
w_1, u(w_2) \rangle.
\]

If $W \cap W_\xi$ is a space $F$ of dimension $3$, then $T_{W_\xi}
(\sigma_3(W)) = \{ u \in \Lambda^2 W_\xi^* \mid u \rest {F \times F}
\equiv 0\}$ according to Lemma \ref{tangent}.  I.e.\ $T_{W_\xi}
(\sigma_3(W))$ is the kernel of the natural map $\Lambda^2 W_\xi^* \to
\Lambda^2 F^*$.

Now let $V_\xi := V/\langle \xi \rangle$.  Then there is a natural
isomorphism $T_\xi \mathbb P^7 \cong V_\xi$.  There is also a natural
isomorphism $W_\xi = \xi \wedge \Lambda^3 V \cong \Lambda^3 V_\xi$.  A
vector $\theta \in T_\xi \mathbb P^7 = V_\xi$ corresponds to a
first-order infinitesimal deformation $\xi + \varepsilon \theta$ of
$\xi$, and hence to the first-order infinitesimal deformation $(\xi +
\varepsilon \theta) \wedge \Lambda^3 V$ of $W_\xi$.  The corresponding
alternating bilinear form on $W_\xi = \xi \wedge \Lambda^3 V$ is
computed by
\[
\langle \xi \wedge \alpha, (\xi+\varepsilon \theta) \wedge \beta
\rangle = \varepsilon (\xi \wedge \alpha \wedge \theta \wedge \beta)
\in \varepsilon (\Lambda^8 V).
\]
If we use the isomorphism $W_\xi \cong \Lambda^3 V_\xi$, then the
corresponding alternating bilinear form on $\Lambda^3 V_\xi$ is
therefore
\[
g_\theta (\alpha,\beta) = - \alpha \wedge \beta \wedge \theta \in
\Lambda^7 V_\xi.
\]

Hence if we think of $F := W \cap W_\xi$ as a three-dimensional
subspace of $\Lambda^3 V_\xi \cong W_\xi$, then $T_{W_\xi}(\mathbb P^7
\cap \sigma_3(W)) \subset V_\xi$ is the orthogonal complement of the
image of $\Lambda^2 F$ in $\Lambda^6 V_\xi \cong V_\xi^*$.  As a
result $(\xi,F,W) \in Y_2$ if and only if the map $\Lambda^2 F \to
\Lambda^6 V_\xi$ is not injective.

Every vector in $\Lambda^2 F$ is of the form $\alpha \wedge \beta$.
We therefore look at the set
\[
T_1(\xi) := \{ (\alpha,\beta) \in (\Lambda^3 V_\xi -\{0\})^2 \mid
\alpha \wedge \beta = 0\}.
\]
We stratify this locus according to the rank of the map $m_\alpha :
V_\xi \to \Lambda^4 V_\xi$ defined by multiplication by $\alpha$.
Generically this map is injective, and the dual multiplication map
$m_\alpha^* : \Lambda^3 V_\xi \to \Lambda^6 V_\xi$ is surjective.  For
$(\alpha,\beta)$ to be in $T_1(\xi)$, $\beta$ must be in the
codimension-$7$ subspace $\ker(m_\alpha^*) \subset \Lambda^3 V_\xi$.
Hence $T_1(\xi)$ contains a stratum of dimension $2n-7$ where
$m_\alpha$ is injective.

The locus $T_1(\xi)$ contains two other strata, but they are of
smaller dimension.  In one of them the $\alpha$ have the form $\alpha'
\wedge \alpha''$ with $\alpha' \in V_\xi$ and $\alpha'' \in \Lambda^2
V_\xi$ indecomposable.  Since $\alpha'$ is really only well-defined up
to a scalar multiple, and $\alpha''$ is really well-defined only in
$\Lambda^2 (V_\xi/\langle \alpha' \rangle)$, these form a locus of
dimension $21 = n-14$ in $\Lambda^3 V$.  For these $\alpha$ the rank
of $m_\alpha$ is $6$, so this stratum in $T_1(\xi)$ has dimension
$2n-20$.

In the final stratum the $\alpha$ are of the form $\alpha_1 \wedge
\alpha_2 \wedge \alpha_3$, a locus of dimension $13 = n-22$.  For
these $\alpha$ the rank of $m_\alpha$ is $4$, so this stratum of
$T_1(\xi)$ has dimension $2n-26$.

Thus $\dim T_1(\xi) = 2n-7$.  As a result the locus
\[
T_2(\xi) := \{ F_0 \in \Gr(2,W_\xi) \mid F_0 = \Span(\alpha,\beta)
\text{ with } \alpha \wedge \beta = 0\}
\]
is of dimension $2n-11$.  Consequently the locus
\[
T_3(\xi) := \{ F \in \Gr(3,W_\xi) \mid F \supset F_0 \text{ for some
$F_0$ in $T_2(\xi)$} \}
\]
is of dimension at most $\dim T_2(\xi)+n-3 = 3n-14$.

For each $F \in T_3(\xi)$, to give a $W \in \OG{n}$ such that $F
\subset W \cap W_\xi$ is equivalent to giving a member of one of the
families of Lagrangian subspaces of $F^\perp/F$.  Consequently, such
$W$ form a family of dimension $(n-3)(n-4)/2$.  As a result, the locus
\[
T_4(\xi) := \{ (F,W) \in \Gr(3,W_\xi) \times \OG{n} \mid F \subset W \cap
W_\xi \text{ and } F \in T_3(\xi)\}
\]
is of dimension at most $(3n-14)+(n-3)(n-4)/2 = n(n-1)/2-8$.  Finally
\[
T_5 := \{ (\xi,F,W) \mid (F,W) \in T_4(\xi) \}
\]
is of dimension at most $n(n-1)/2-1$.  But the part of $T_5$ where $F
= W \cap W_\xi$ is exactly the part of $Y_2$ lying in $Y-Y_1$.  So
$\dim Y_2 < n(n-1)/2 = \dim (\OG{n})$.

As a result $Y_2$ does not dominate $\OG{n}$, and $\Upsilon \subset
\OG{n}$ is a nonempty open subset.  As before, $\Upsilon$ must then
contain a rational point because the base field is infinite, and we
conclude.
\end{proof}

%

\section{Extension of scalars and non-Pfaffian examples}
\label{sect.real}

In this section we prove the following extension of Theorem
\ref{counterexample}.  It will allow us to construct subcanonical
subschemes of codimension $3$ in $\mathbb P^N_k$ which are not
Pfaffian over $k$ but are Pfaffian over a finite extension of $k$.

\begin{theorem}
\label{sym.bilin.thm}
Let $k$ be an infinite field, let $E$ be a $k$-vector space of
dimension $d \geq 1$, and let $b : E \times E \to k$ be a
nondegenerate symmetric bilinear form.  There exists a nonsingular
subcanonical fourfold $X \subset \mathbb P^7_k$ of degree $(1000d^3 +
8d)/3$ with $K_X = (20d-8)H$ such that the cup product pairing
\eqref{cup} on the middle cohomology
\[
H^2(\mathcal O_X(10d-4)) \times H^2(\mathcal O_X(10d-4)) \to
H^4(\mathcal O_X(20d-8)) \cong k
\]
may be identified with $b : E \times E \to k$.
\end{theorem}

\begin{proof}
If $e_1, \dots , e_d$ is a basis of $E$ and $b_{ij} := b(e_i,e_j)$,
then the quadratic form on $E \otimes \Lambda^4 V$ given by
\[
Q \left( \sum_i e_i \otimes u_i \right) := \sum_i b_{ii} \,
u_i^{(2)} + \sum_{i<j} b_{ij} \, u_i \wedge u_j
\]
is nondegenerate and has the property that if $L$ is a Lagrangian
subspace of $(\Lambda^4 V,q)$, then $E \otimes L$ is a Lagrangian
subspace of $(E \otimes \Lambda^4 V, Q)$.  Consequently $E \otimes
\Omega^4_{\mathbb P^7}(4)$ is a Lagrangian subbundle of $(E \otimes
\Lambda^4 V \otimes \mathcal O_{\mathbb P^7},Q)$.  Also, if we make a
standard identification $\Lambda^4 V^* \cong \Lambda^4 V$, then the
symmetric bilinear form associated to $Q$ is $b \otimes 1: E
\otimes \Lambda^4 V \to E^* \otimes \Lambda^4 V$.

We now choose a general Lagrangian subspace $W \subset (E \otimes
\Lambda^4 V,Q)$ and apply the construction of Theorem \ref{split}.
This produces a subcanonical fourfold $X \subset \mathbb P^7_k$ with
$K_X = (20d-8)H$.  The degree of $X$ can be computed from the formula
after Theorem \ref{split}.  If $\Char(k) =0$, then $X$ is nonsingular
by Theorem \ref{split}.  If $k$ is infinite of positive
characteristic, $X$ is nonsingular by an argument similar to Lemma
\ref{transversal} whose details we leave to the reader.

Let $N := 10d-4$.  The isomorphism $\eta : \mathcal O_X(N) \to
\omega_X(-N)$ induces an isomorphism
\begin{equation}
\label{serre.dual}
H^2(\mathcal O_X(N)) \xrightarrow{\sim} H^2(\omega_X(-N)) \cong
H^2(\mathcal O_X(N))^*
\end{equation}
which corresponds to the cup product pairing of Serre duality
\eqref{cup}.  Our problem is to identify this map.

Now twisting the diagram \eqref{dual.diag} of Theorem \ref{split}
produces a commutative diagram with exact rows (where we write
$\mathcal O := \mathcal O_{\mathbb P^7}$ to keep things within the
margins):
\begin{small}
\begin{diagram}[LaTeXeqno]
\label{two.lines}
0 & \rTo & \mathcal O(-N-8) & \rTo^f & E \otimes \Omega^4_{\mathbb
P^7} & \rTo^\psi & W^* \otimes \mathcal O(-4) & \rTo^g & \mathcal O(N)
& \rOnto & \mathcal O_{X}(N) \\
&& \dSame && \dTo <\phi && \dTo >{\phi^*}  && \dSame && \dSame \\
0 & \rTo & \mathcal O(-N-8) & \rTo^{g^*} & W \otimes \mathcal O(-4) &
\rTo^{-\psi^*} & E^* \otimes \Omega^3_{\mathbb P^7} & \rTo^{f^*} &
\mathcal O(N) & \rOnto & \omega_{X}(-N)
\end{diagram}%
\end{small}%
The rows of the diagram produce natural isomorphisms $H^2(\mathcal
O_X(N)) \cong E$ and $H^2(\omega_X(-N)) \cong E^*$.  Moreover, the
Serre duality pairing
\[
H^2(\mathcal O_X(N)) \times H^2(\omega_X(-N)) \to H^4(\omega_X)
\xrightarrow{\tr} k
\]
corresponds to the canonical pairing $E \times E^* \to k$ since the
bottom resolution is dual to the top resolution, so the
hypercohomology of the bottom resolution is naturally Serre dual to
the hypercohomology of the top resolution.  It remains to identify the
map $H^2(\mathcal O_X(N)) \to H^2(\omega_X(-N))$ with the map $b : E
\to E^*$.

We do this by using the following commutative diagram with exact
rows. 
\begin{footnotesize}
\begin{diagram}[LaTeXeqno]
\label{three.lines}
0 & \rTo & \mathcal O(-N-8) & \rTo^f & E \otimes \Omega^4_{\mathbb
P^7} & \rTo^\psi & W^* \otimes \mathcal O(-4) & \rTo^g & \mathcal O(N)
& \rOnto & \mathcal O_{X}(N) \\
&& \dSame && \dTo <\phi && \dTo >{\alpha}  && \dSame && \dSame \\
0 & \rTo & \mathcal O(-N-8) & \rTo & W \otimes \mathcal O(-4) &
\rTo^{-\beta} & E \otimes \Omega^3_{\mathbb P^7} & \rTo &
\mathcal O(N) & \rOnto & \mathcal O_X(N) \\
&& \dSame && \dSame && \dTo >{b \otimes 1} && \dSame && \dTo >\eta
<\cong \\
0 & \rTo & \mathcal O(-N-8) & \rTo^{g^*} & W \otimes \mathcal O(-4) &
\rTo^{-\psi^*} & E^* \otimes \Omega^3_{\mathbb P^7} & \rTo^{f^*} &
\mathcal O(N) & \rOnto & \omega_{X}(-N)
\end{diagram}%
\end{footnotesize}%
The maps $\alpha$ and $\beta$ in this diagram come from the
identification $E \otimes \Lambda^4 V \cong W^* \oplus W$, which
yields a short exact sequence
\begin{small}
\begin{equation}
\label{WW.ex.seq}
0 \to E \otimes \Omega_{\mathbb P^7}^4 \xrightarrow{\sm{\psi\\ \phi}}
\left( W^* \otimes \mathcal O (-4) \right) \oplus \left( W \otimes
\mathcal O (-4) \right) \xrightarrow{\sm{\alpha& \beta}} E \otimes
\Omega_{\mathbb P^7}^3 \to 0.
\end{equation}%
\end{small}%
The exactness of \eqref{WW.ex.seq} implies that the top middle square
of diagram \eqref{three.lines} commutes, and that $\phi$ induces an
isomorphism $\ker(\psi) \xrightarrow{\sim} \ker(\beta)$, while
$\alpha$ induces an isomorphism $\coker(\psi) \xrightarrow{\sim}
\coker(\beta)$.  As a result the second row of the diagram is exact,
and the squares in the top row are all commutative.

To see that the squares in the bottom row are all commutative, recall
how the maps $\alpha$, $\beta$, $\phi^*$ and $\psi^*$ are defined.
The inclusions of $W^*$ and $W$ in $W^* \oplus W \cong E \otimes
\Lambda^4 V$, together with the bilinear form and the Koszul complex,
give rise to a diagram with a commutative square
\begin{diagram}[loose,LaTeXeqno,PS]
\label{with.square}
W \otimes \mathcal O(-4) & \rInto^{i} & E \otimes \Lambda^4 V \otimes
\mathcal O(-4) & \rTo^{b \otimes \boldsymbol 1} & E^* \otimes
\Lambda^4 V \otimes \mathcal O(-4) \\
& \ruInto >{i'} & \dTo >\pi && \dTo >{\pi'} \\
W^* \otimes \mathcal O(-4) && E \otimes \Omega_{\mathbb P^7}^3 & \rTo
^{b \otimes 1_\Omega} & E^* \otimes \Omega_{\mathbb P^7}^3
\end{diagram}%
such that 
\begin{align*}
\alpha & = \pi i', & \beta & = \pi i, & \psi^* & = \pi' (b \otimes
\boldsymbol 1) i, & \phi^* & = \pi' (b \otimes \boldsymbol 1) i'.
\end{align*}
Since the square in \eqref{with.square} is commutative, we get $\psi^*
= (b\otimes 1_\Omega)\beta$ and $\phi^* = (b \otimes 1_\Omega)\alpha$.
The first equality means that the bottom set of squares in diagram
\eqref{three.lines} commutes.  The second equality means that the
composite of the two chain maps of diagram \eqref{three.lines} is the
chain map of diagram \eqref{two.lines}.

Now the first and second rows of diagram \eqref{three.lines} produce
compatible identifications $H^2(\mathcal O_X) \cong E$.  And the
second and third rows of diagram \eqref{three.lines} show that the map
$H^2(\mathcal O_X(N)) \to H^2(\omega_X(-N))$ induced by $\eta$ may be
identified with $b : E \to E^*$.  This completes the proof of the
theorem.
\end{proof}

As an example, we may apply Theorem \ref{sym.bilin.thm} with $b$ a
positive definite symmetric bilinear form on $\mathbb R^2$.  There are
no Lagrangian subspaces of $(\mathbb R^2,b)$, but there are Lagrangian
subspaces of $(\mathbb C^2,b_{\mathbb C})$.  So by the Pfaffianness
criterion of Theorem \ref{criterion}, we get the following corollary.

\begin{corollary}
\label{real.counterexample}
There exists a nonsingular subcanonical fourfold $Y_{\mathbb R}
\subset \mathbb P^7_{\mathbb R}$ of degree $2672$ with $K_Y = 32 H$
which is not Pfaffian over $\mathbb R$ but whose complexification
$Y_{\mathbb C} \subset \mathbb P^7_{\mathbb C}$ is Pfaffian over
$\mathbb C$.
\end{corollary}

\section{Reisner's example in $\mathbb P^5$}
\label{sect.rp2}

Another example of a codimension $3$ subcanonical subscheme $Z\subset
\mathbb P^5$ is provided by the Stanley-Reisner ring of the minimal
triangulation of the real projective plane $\mathbb R\mathbb P^2$ in
characteristic $2$.  This is Reisner's original example of a monomial
ideal whose minimal free resolution depends on the characteristic of
the base field.  See Hochster \cite{Hochster}, Reisner \cite{Reisner},
Stanley \cite{Stanley}.

Recall that if $\Delta$ is a simplicial complex on the vertex set $E=
\{e_0,\ldots e_n\}$, the corresponding face variety $X(\Delta) \subset
\mathbb P^n_k$ in the sense of Stanley, Hochster and Reisner (see
\cite{Stanley} for more details) is the locus defined by the ideal
$I_{\Delta}$ generated by monomials corresponding to the simplexes
(faces) not contained in $\Delta$.  Geometrically, $X(\Delta)$ is a
``projective linear realization'' of $\Delta$.  It is the union of
linear subspaces of $\mathbb P^n_k$ corresponding to the simplices in
$\Delta$, where for each $m$-simplex $\{ e_{i_0}, \dots, e_{i_m} \}$
one includes the $m$-plane spanned by the vertices $P_{i_0}, \dots,
P_{i_m}$ of the standard coordinate system on $\mathbb P^n_k$.

If the topological realization $|\Delta|$ is a manifold, then
$X(\Delta)$ is a locally Gorenstein scheme with $\omega_{X(\Delta)}
^{\otimes 2} \cong \mathcal O_{X(\Delta)}$.  Moreover,
$\omega_{X(\Delta)}$ is trivial if and only if $|\Delta|$ is
orientable over the field $k$.  The cohomology of $X(\Delta)$ is given
by the formula $H^i(\mathcal O_{X(\Delta)}) = H^i(\Delta, k)$.  The
Hilbert polynomial of $X(\Delta)$, which coincides with the Hilbert
function for strictly positive values, is completely determined by the
combinatorial data (see \cite{Stanley} for more details).


\begin{figure}[hbtp]
\begin{center}
\leavevmode
\epsfbox{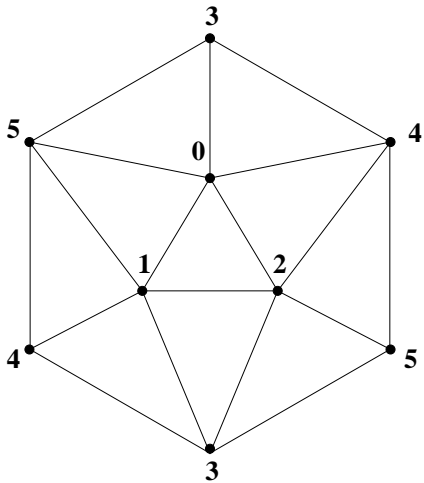}
\caption{The Minimal Triangulation of $\mathbb {RP}^2$}
\end{center}
\end{figure}

We apply this theory with $\Delta$ the triangulation of $\mathbb
{RP}^2$ given in Figure 1, and with $k$ a field of characteristic $2$.
Then $X := X(\Delta)$ is a union of ten $2$-planes in $\mathbb P^5_k$
corresponding to the ten triangles.  Moreover, $\omega_X \cong
\mathcal O_X$ since $\mathbb {RP}^2$ is orientable in characteristic
$2$.  The cohomology is given by
\begin{align*}
h^0(\mathcal O_X) & = 1, & h^1(\mathcal O_X) & = 1, & h^2(\mathcal
O_X) & = 1,
\end{align*}
since this is the cohomology of $\mathbb {RP}^2$ over a field of
characteristic $2$.  We may now prove

\begin{proposition}[$\Char\ 2$]
\label{RP2}
Let $X \subset \mathbb P^5$ be the face variety corresponding to the
triangulation of $\mathbb {RP}^2$ given in Figure \textup 1.  Then $X$
is subcanonical of codimension $3$ but not Pfaffian.
\end{proposition}

\begin{proof}
The subvariety $X$ is a locally Gorenstein surface with $\omega_X =
\mathcal O_X$ but $h^1(\mathcal O_X) = 1$.  Thus it is subcanonical
but not Pfaffian by Theorem \ref{criterion}.  
\end{proof}

In addition, $X \subset \mathbb P^5$ is linearly normal (all
$X(\Delta)$ are) and is not contained in a hyperquadric (since all
pairs of vertices are joined by an edge in $\Delta$).  This will allow
us to conclude below that $X$ is a degenerate member of a family of
nonclassical Enriques surfaces (see Proposition \ref{converse.P5} and
\eqref{ext.monoms}).  Over the complex numbers this is a Type III
degeneration in Morrison \cite{Mor} Corollary 6.2 (iii a).  See also
Symms \cite{Symms} and Altmann-Christophersen (in preparation) for the
deformation theory of this scheme.

\section{Enriques surfaces in $\mathbb P^5$}
\label{sect.Enriques}

In their extension to characteristic $p$ of the
Enriques-Castelnuovo-Kodaira classification of surfaces not of general
type, Bombieri and Mumford \cite{BM} characterize Enriques surfaces as
nonsingular projective surfaces $X$ with numerically trivial canonical
class $K_X \equiv 0$ and with $\chi(\mathcal O_X) = 1$.  In
characteristic $2$, Enriques surfaces are divided into three types:

\begin{center}
\begin{small}
\begin{tabular}{|c||c|c|c|c|} \hline
& $h^1(\mathcal O_X)$ & $h^2(\mathcal O_X)$ &
canonical class & $F : H^1(\mathcal O_X)\to H^1(\mathcal O_X)$
\\ \hline\hline
classical & $0$ & $0$ & $2K_X = 0$, $K_X \neq 0$ &  --- \\
\hline 
$\mmu$ & $1$ & $1$ & $K_X = 0$ & injective \\ \hline
$\aalpha$ & $1$ & $1$ & $K_X = 0$ & $0$ \\ \hline
\end{tabular}%
\end{small}%
\end{center}%
The map $F: H^1(\mathcal O_X) \to H^1(\mathcal O_X)$ is the action of
Frobenius on the cohomology.  There is a voluminous literature on
Enriques surfaces which we do not cite here.  Our Pfaffianness
criterion Theorem \ref{criterion} yields immediately

\begin{proposition}[$\Char \ 2$]
If $X \subset \mathbb P^5$ is a nonclassical Enriques surface, then
$X$ is subcanonical but not Pfaffian.
\end{proposition}

A {\em Fano polarization} on an Enriques surface $X$ is a numerically
effective divisor $H$ such that $H^2 = 10$ and such that $H \cdot F
\geq 3$ for every numerically effective divisor $F$ with $F^2 = 0$.
It is well known that every Enriques surface admits a Fano
polarization, and that for any Fano polarization $H$ the linear system
$\lvert H \rvert$ is base-point-free and defines a map $\varphi_H : X
\to \mathbb P^5$ which is birational onto its image.  The image
$\varphi_H(X) \subset \mathbb P^5$ is called a {\em Fano model} of
$X$.  A surface $Y \subset \mathbb P^5$ is a Fano model of an Enriques
surface if and only if it is of degree $10$ with at worst rational
double points as singularities and has $\chi(\mathcal O_Y) = 1$ and
numerically trivial canonical class $K_Y \equiv 0$.

An Enriques surface is called {\em nodal} if it contains a smooth
rational curve $C$ with $C^2 = -2$.  Fano models of unnodal Enriques
surfaces are necessarily smooth.  The following theorem was proven in
characteristic $\neq 2$ by Cossec \cite{Cossec}.  Igor Dolgachev has
informed us that the theorem also holds in characteristic $2$ using
arguments similar to Dolgachev-Reider \cite{DR}.

\begin{theorem}
\label{unnodal}
Let $H$ be a Fano polarization of an Enriques surface $X$.  Then the
following are equivalent\textup: \textup{(a)} $X$ is nodal, and
\textup{(b)} one of the Fano models $\varphi_H(X)$ or
$\varphi_{H+K_X}(X)$ is contained in a hyperquadric.
\end{theorem}

We will now give an algebraic method for realizing Fano models of
unnodal nonclassical Enriques surfaces as degeneracy loci associated
to Lagrangian subbundles of an orthogonal bundle on $\mathbb P^5$.  We
will need the following technical lemmas about divided squares.  Both
lemmas hold for even $m$ in all characteristics.


\begin{lemma}
\label{div.square.lemma}
Suppose that $R$ is a commutative algebra over a field of
characteristic $2$, that $m \geq 3$ is an odd integer, and that $E$ a
free $R$-module of rank $2m$.  Then there exists a well-defined
divided square operation $\Lambda^m E \to \Lambda^{2m} E \cong R$
whose formula, with respect to any basis $e_1,\dots, e_m$, is given by
\begin{equation}
\label{div.square}
\left( \sum x_{i_1 \dots i_m} e_{i_1 \dots i_m} \right)^{(2)} :=
\sum_{\substack {i_1<\dotsb < i_m \\ j_1 < \dotsb < j_m \\ i_1 < j_1
\\ \{i_k\} \cap \{j_\ell\} = \varnothing}} x_{i_1 \dotsb i_m} x_{j_1
\dotsb j_m} e.
\end{equation}
\end{lemma}

\begin{proof}
It is enough to prove the lemma for $R$ a field, since one can reduce
to the case where the coefficients $x_{i_1 \dots i_m}$ are independent
indeterminates, and $R$ is the field of rational functions in these
indeterminates.

We have to show that \eqref{div.square} is invariant under changes of
basis in $\GL(E) = \GL_n(R)$.  But since $R$ is now assumed to be a
field, $\GL_n(R)$ is generated by three special kinds of changes of
basis: permutations of the basis elements, operations of the form $e_i
\rightsquigarrow \lambda e_i$, and operations of the form $e_j
\rightsquigarrow e_j + \alpha e_i$.  The formula \eqref{div.square} is
almost trivially invariant under the first two kinds of operations,
and it is invariant under the third kind of operation because the
portions
\[
x_{ij k_1 \dots k_{m-2}} x_{\ell_1 \dots \ell_m} \qquad \text{and}
\qquad  x_{i k_1 \dots k_{m-1}} x_{j \ell_1 \dots \ell_{m-1}} + x_{i
\ell_1 \dots \ell_{m-1}} x_{j k_1 \dots k_{m-1}}
\]
do not change.
\end{proof}

\begin{lemma}
\label{unique.div.square}
Suppose that $m \geq 3$ is an odd integer, that $V$ is a vector space
of dimension $2m$ over a field of characteristic $2$, and that $Q$ is
a nonzero quadratic form on $\Lambda^m V$.  Then $Q(x \wedge w) = 0$
for all $x \in V$ and $w \in \Lambda^{m-1}V$ if and only if $Q$ is a
constant multiple of the divided-square quadratic form of Lemma
\ref{div.square.lemma}.
\end{lemma}

\begin{proof}
Suppose first that $Q$ is a multiple of the divided-square quadratic
form, and that $0 \neq x \in V$ and $w \in \Lambda^{m-1}V$.   Choosing
a basis of $V$ of the form $e_1 = x, e_2, \dots, e_m$ and applying
formula \eqref{div.square}, one sees easily that $Q(x\wedge w) = 0$.

Conversely, suppose that $Q(x\wedge w) = 0$ for all $x \in V$ and $w
\in \Lambda^{m-1}V$.  Write $Q(\sum y_I e_I) := \sum_{I \preceq J}
\alpha_{IJ} y_I y_J$.  The hypothesis implies that $Q(e_I) =
\alpha_{II} = 0$ for all $I$.  It also implies that each $Q(e_{iA} +
e_{iB}) = \alpha_{iA,iB} = 0$.  Hence $\alpha_{IJ} = 0$ if the
multi-indices $I$ and $J$ are not disjoint.  Finally
$Q((e_i+e_j)\wedge (e_A + e_B)) = \alpha_{iA,jB} + \alpha_{iB,jA} =
0$.  In other words $\alpha_{IJ} = \alpha_{KL}$ if $I$ and $J$ are
disjoint, and $K$ and $L$ are disjoint, and $I$ and $K$ contain
exactly $m-1$ common indices.  One then easily deduces that all the
$\alpha_{IJ}$ with $I$ and $J$ disjoint must be equal.  But then $Q$
is the same as the quadratic form of \eqref{div.square} up to a
constant.
\end{proof}

Now let $V = H^0(\mathcal O_{\mathbb P^5}(1))^*$, and consider the
exact sequence of vector bundles on $\mathbb P^5$.
\begin{equation}
\label{trunc.Koszul}
0 \to \Omega_{\mathbb P^5}^3(3) \to \Lambda^3 V \otimes \mathcal
O_{\mathbb P^5} \to \Omega_{\mathbb P^5}^2(3) \to 0.
\end{equation}
As in \eqref{Euler}, if $\xi \in V$, then the fiber of
$\Omega_{\mathbb P^5}^3(3)$ over $\overline \xi \in \mathbb P^5 =
\mathbb P(V)$ may be identified with the subspace $\xi \wedge
\Lambda^2 V \subset \Lambda^3 V$.  Therefore Lemmas
\ref{div.square.lemma} and \ref{unique.div.square} have the following
corollary.

\begin{corollary}[$\Char \ 2$]
\label{unique.quadratic}
There exists a nondegenerate quadratic form $Q$ on $\Lambda^3 V$,
unique up to a constant multiple, such that $\Omega_{\mathbb
P^5}^3(3)$ is a Lagrangian subbundle of the orthogonal bundle
$(\Lambda^3 V \otimes \mathcal O_{\mathbb P^5},Q)$.
\end{corollary}

We now apply the construction of Theorem \ref{split} to the Lagrangian
subbundle $\Omega_{\mathbb P^5}^3(3) \subset \Lambda^3 V \otimes
\mathcal O_{\mathbb P^5}$.  We get the following result.

\begin{theorem}[$\Char\ 2$]
\label{Fano.existence}
If $W \subset (\Lambda^3 V,Q)$ is a general Lagrangian subspace in the
family opposite to the one containing the fibers of $\Omega^3_{\mathbb
P^5}(3)$, then the degeneracy locus
\begin{equation}
\label{ZW}
Z_W := \{x \in \mathbb P^5 \mid \dim \left[\Omega^3_{\mathbb
P^5}(3)(x) \cap W \right] \geq 3 \}
\end{equation}
is a smooth Fano model of a nonclassical unnodal Enriques surface with 
symmetrically  quasi-isomorphic resolutions
\begin{small}
\begin{diagram}[LaTeXeqno]
\label{P5.diag}
0 & \rTo & \mathcal O_{\mathbb P^5}(-6) & \rTo & \Omega^3_{\mathbb
P^5} & \rTo^\psi & W^* \otimes \mathcal O_{\mathbb P^5}(-3) & \rTo &
\mathcal O_{\mathbb P^5} & \rOnto & \mathcal O_{Z_W} \\
&& \dSame && \dTo <\phi && \dTo >{\phi^*} && \dSame && \dTo >\eta
<\cong \\
0 & \rTo & \mathcal O_{\mathbb P^5}(-6) & \rTo & W \otimes \mathcal
O_{\mathbb P^5}(-3) & \rTo^{-\psi^*} & \Omega^2_{\mathbb P^5} & \rTo &
\mathcal O_{\mathbb P^5} & \rOnto & \omega_{Z_W}
\end{diagram}
\end{small}%
\end{theorem}

\begin{proof}
The general $Z_W$ is a subcanonical surface of degree $10$ with
$\omega_{Z_W} \cong \mathcal O_{Z_W}$ and with the resolutions
\eqref{P5.diag} according to Theorem \ref{split}.  It is smooth by an
argument similar to Lemma \ref{transversal}.  Since $K_{Z_W}= 0$ and
$h^i(\mathcal O_{Z_W}) = 1$ for $i=0,1,2$, it is a nonclassical
Enriques surface by the Bombieri-Mumford classification.  Since the
degree is $10$, it is a Fano model.  Moreover, according to the
resolution, it is not contained in a hyperquadric.  Hence it is
unnodal by Theorem \ref{unnodal}.
\end{proof}

The theorem has the following converse.

\begin{proposition}[$\Char\ 2$]
\label{converse.P5}
Let $X \subset \mathbb P^5$ be a locally Gorenstein subscheme of
dimension $2$ and degree $10$, with $\omega_X \cong \mathcal O_X$,
with $h^i(\mathcal O_X) = 1$ for $i=0,1,2$, and which is linearly
normal and not contained in a hyperquadric.  Let $V = H^0(\mathcal
O_{\mathbb P^5}(1))^*$, and let $Q$ be the divided square on
$\Lambda^3 V$.  Then there exists a unique Lagrangian subspace $W
\subset (\Lambda^3 V,Q)$ such that $X = Z_W$ \textup(cf.\ Theorem
\ref{Fano.existence}\textup) with locally free resolutions as in
\eqref{P5.diag}.
\end{proposition}

For example, the face variety $X = X(\Delta)$ of Proposition \ref{RP2}
satisfies all the hypotheses of the proposition and is defined as $X =
Z_W$ for the Lagrangian subspace $W \subset \Lambda^3 V$ spanned by
the exterior monomials
\begin{equation}
\label{ext.monoms}
e_{013},\ e_{014},\ e_{023},\ e_{025},\ e_{045},\ e_{124},\ e_{125},\ 
e_{135},\ e_{234},\ e_{345}.
\end{equation}
%

Any $Z_W$ in \eqref{ZW} which is of dimension $2$ satisfies all the
hypotheses of Proposition \ref{converse.P5}.

The following refinement of the Castelnuovo-Mumford Lemma is needed
for the proof of Proposition \ref{converse.P5}.  Its proof is left to
the reader.

\begin{lemma}
\label{refined.CM}
Let $\mathcal M$ be a coherent sheaf on $\mathbb P^n$ with no nonzero
skyscraper subsheaves, and let $0 \leq q \leq n-1$ and $m$ be
integers.  Suppose that $\mathcal M$ is $(m+1)$-regular, that
$H^i(\mathcal M(m-i)) = 0$ for $i \geq n-q$, and that the
comultiplication maps
\[
H^i(\mathcal M(m-1-i)) \to V \otimes H^i(\mathcal M(m-i))
\]
are surjective for $1 \leq i < n-q$.  Then for any $0 \leq j \leq q$,
the module of $j$-th syzygies of $H^0_*(\mathcal M)$ is generated in
degrees $\leq m+j$.
\end{lemma}

\begin{proof}[Proof of Proposition \ref{converse.P5}]
The Hilbert polynomial of $X$ is $\chi(\mathcal O_X(t)) = 5t^2+1$ by
Riemann-Roch.  The reader may easily verify that the sheaves $\mathcal
O_X(t)$ have cohomology as given in the following table:

\begin{center}
\begin{small}
\begin{tabular}{|c||c|c|c|} \hline
& $h^0(\mathcal O_X(t))$ & $h^1(\mathcal O_X(t))$ & $h^2(\mathcal
O_X(t))$ \\ \hline\hline
$t > 0$ & $5t^2+1$ & $0$ & $0$ \\
\hline 
$t = 0$ & $1$ & $1$ & $1$ \\ \hline
$t < 0$ & $0$ & $0$ & $5t^2+1$ \\ \hline
\end{tabular}%
\end{small}%
\end{center}%
%

Now according to the table $\mathcal O_X$ is $3$-regular; the only
group obstructing its $2$-regularity is $H^2(\mathcal O_X)$; and the
comultiplication $H^2(\mathcal O_X(-1)) \to V \otimes H^2(\mathcal
O_X)$ is dual to the multiplication $V^* \otimes H^0(\mathcal O_X) \to
H^0(\mathcal O_X(1))$ and thus is an isomorphism.  So we may apply
Lemma \ref{refined.CM} to see that all the generators of
$H^0_*(\mathcal O_X)$ are in degrees $\leq 2$, and all the relations
are in degrees $\leq 3$.  But for all $t \leq 2$ the restriction maps
$H^0(\mathcal O_{\mathbb P^5}(t)) \to H^0(\mathcal O_X(t))$ are
bijections.  So the generators and relations of $H_*^0(\mathcal
O_{\mathbb P^5})$ and of $H^0_*(\mathcal O_X)$ in degrees $\leq 2$ are
the same---one generator in degree $0$ and no relations.  Hence
$H^0_*(\mathcal O_X)$ has only a single generator in degree $0$, and
its relations are all in degree $3$.  In other words $H^1_*(\mathcal
I_X) = 0$, and the homogeneous ideal $I(X)$ is generated by a
$10$-dimensional vector space $W^*$ of cubic equations.

We now go through the argument in \cite{EPW} Theorem 6.2 which
constructs a pair of symmetrically quasi-isomorphic locally free
resolutions of a subcanonical subscheme of codimension $3$.  We have
an exact sequence
\[
0 \to \mathcal K \to W^* \otimes \mathcal O_{\mathbb P^5}(-3) \to
\mathcal O_{\mathbb P^5} \to \mathcal O_X \to 0.
\]
There is a natural isomorphism 
\[
\Ext^1_{\mathcal O_{\mathbb P^5}}(\mathcal K, \mathcal O_{\mathbb
P^5}(-6)) \cong H^4(\mathcal K)^* \cong H^2(\mathcal O_X)^* \cong k,
\]
and a nonzero member of this group induces an extension which we can
attach to give a locally free resolution
\[
0 \to \mathcal O_{\mathbb P^5}(-6) \to \mathcal F \to W^* \otimes
\mathcal O_{\mathbb P^5}(-3) \to \mathcal O_{\mathbb P^5} \to \mathcal
O_X \to 0.
\]
The intermediate cohomology of $\mathcal F$ comes from the
intermediate cohomology of $X$ and is given by $H^i_*(\mathcal F) = 0$
for $i=1,2,4$ and $H^3_*(\mathcal F) \cong H^1_*(\mathcal O_X) = k$.
According to Horrocks' Theorem \cite{OSS}, this means that $\mathcal F
\cong \Omega^3_{\mathbb P^5} \oplus \bigoplus \mathcal O_{\mathbb
P^5}(a_i)$.  Since $\mathcal F$ and $\Omega^3_{\mathbb P^5}$ are both
of rank $10$, they are isomorphic.  This gives the top locally free
resolution of \eqref{P5.diag}.  The rest of diagram \eqref{P5.diag} is
now constructed as in the proof of \cite{EPW} Theorem 6.2.

It remains to identify the quadratic spaces $W^* \oplus W$ and
$\Lambda^3 V$, and to show that $W$ is uniquely determined by $X$.
However, \eqref{P5.diag} gives us an exact sequence
\[
0 \to \Omega_{\mathbb P^5}^3(3) \xrightarrow {\sm{\psi \\ \phi}} (W^*
\oplus W) \otimes \mathcal O_{\mathbb P^5} \xrightarrow{\sm{\phi^* &
\psi^*}} \Omega_{\mathbb P^5}^2(3) \to 0.
\]
This is the unique nontrivial extension of $\Omega_{\mathbb P^5}^2(3)$
by $\Omega_{\mathbb P^5}^3(3)$, and thus coincides with the truncated
Koszul complex \eqref{trunc.Koszul}.  There is therefore a natural
identification of $W^* \oplus W$ with $\Lambda^3 V$ which identifies
the hyperbolic quadratic form on $W^* \oplus W$ with a quadratic form
$Q$ on $\Lambda^3 V$ for which $\Omega_{\mathbb P^5}^3(3) \subset
\Lambda^3 V \otimes \mathcal O_{\mathbb P^5}$ is a Lagrangian
subbundle.  By Corollary \ref{unique.quadratic}, $Q$ is necessarily
the divided square quadratic form.   

The identification of the quadratic spaces $W^* \oplus W$ and
$\Lambda^3 V$ is unique up to homothety, once the morphisms $\psi,
\phi$ are chosen.  However, the diagram \eqref{P5.diag} is unique up
to an alternating homotopy $\alpha : W^*\otimes \mathcal O_{\mathbb
  P^5} \to W \otimes \mathcal O_{\mathbb P^5}$.  One may check that
this alternating homotopy does not change the Lagrangian subspace $W
\subset \Lambda^3 V$ but only the choice of the Lagrangian complement
$W^*$ (cf.\ \cite{EPW} \S 5).  Thus $W$ is uniquely determined by $X$.
This completes the proof.
\end{proof}

\begin{corollary}[$\Char\ 2$]
\label{moduli}
There is a universal family of Fano models in $\mathbb P^5 = \mathbb
P(V)$ of unnodal nonclassical Enriques surfaces parametrized by an
$\SL(V)$-invariant Zariski open subset $U \subset \og_{20}(\Lambda^3
V)$.  The $\SL(V)$-orbits correspond to isomorphism classes of
Fano-polarized unnodal nonclassical Enriques surfaces $(X,H)$.
\end{corollary}

In particular, unnodal nonclassical Enriques surfaces form an
irreducible family.  We will see below that its general member is a
$\mmu$-surface (Proposition \ref{Frob.locus}).

The universal family of unnodal nonclassical Enriques surfaces has the
following universal locally free resolution over $\mathbb G :=
\og_{20}(\Lambda^3 V)$.  Let $\mathcal S \subset \Lambda^3 V \otimes
\mathcal O_{\mathbb G}$ by the universal Lagrangian subbundle on
$\mathbb G$.  According the machinery of \cite{EPW} Theorem 2.1, the
two subbundles $p^* \Omega_{\mathbb P^5}^3(3)$ and $q^* \mathcal S$ on
$\mathbb P^5 \times \mathbb G$ (where $p$ and $q$ are the two
projections) define a degeneracy locus $Z$, which one can verify has
the expected codimension $3$ using an incidence variety argument.
Therefore $Z \subset \mathbb P^5 \times \mathbb G$ has a locally free
resolution
\begin{equation}
\label{universal}
0 \to \mathcal O(-6,-2) \to (p^* \Omega_{\mathbb P^5}^3) (0,-1) \to
(q^*\mathcal S^*)(-3,-1) \to \mathcal O \to \mathcal O_Z \to 0.
\end{equation}

\begin{remark} 
\label{geometric}
(a) Barth and Peters have shown that a generic Enriques surface over
$\mathbb C$ has exactly $2^{14} \cdot 3 \cdot 17 \cdot 31$ distinct
Fano polarizations (\cite{BP} Theorem 3.11).

(b) The fact that the Fano model of a unnodal Enriques surface is
defined by $10$ cubic equations was first observed by F.\ Cossec (cf.\ 
\cite{Cossec} \cite{CD}).

(c) ($\Char\ 2$)
The following geometrical construction of a general Fano-polarized
unnodal $\mmu$-surface is inspired by Mori-Mukai (\cite{MM}
Proposition 3.1).  Let $E$ be an elliptic curve with Hasse invariant
one, and let $i: E \to E$ be the involution given by translating by
the unique nontrivial $2$-torsion point.  Let $\pi : E \to E' :=
E/\langle i \rangle$ be the quotient map.  Let $D'$ be a divisor of
degree 3 on $E'$ and let $D:=\pi^*(D')$. The linear system $|D|$ gives
an embedding $E \hookrightarrow \mathbb P^5$ with image an elliptic
normal sextic curve.  The involution $i$ on $E$ extends to an
involution of $\mathbb P^5$.  As in Mori-Mukai \cite{MM}, (3.1.1) we
may choose three general $i$-invariant hyperquadrics $Q_1$, $Q_2$,
$Q_3$ containing $E$, such that $X := Q_1 \cap Q_2 \cap Q_3 \subset
\mathbb P^5$ is a smooth K3 surface, such that $i \rest{X}$ has no
fixed points.  Let $L$ be the hyperplane divisor class on $X$. Then
$S:= X/\langle i \rangle$ is an Enriques surface of type $\mmu$ (since
it has an \'etale double cover by a K3 surface, cf.\ \cite{BM}).  The
divisor class $L+E$ descends to a Fano polarization $H$ on $S$, and
the elliptic curve $E \subset X$ descends to $E' \subset S$.  By
\cite{CV} or \cite{DR} Theorem 1, (which Igor Dolgachev informs us
also holds in characteristic $2$), $S$ is unnodal if and only if
$|H-2E'|=\emptyset$.  But this follows since $|L-E|=\emptyset$ because
$E \subset \mathbb P^5$ is nondegenerate.  An easy dimension count
shows that this construction gives rises to $10$ moduli of Fano
polarized unnodal $\mmu$ Enriques surfaces, and thus to a general one.
\end{remark}

\section{Calculating the action of Frobenius}
\label{sect.Frob}

We would like to give a method for identifying when a Fano model of a
nonclassical Enriques surface constructed in Theorem
\ref{Fano.existence} is of type $\mmu$ or type $\aalpha$.  This
involves computing explicitly the action of Frobenius on the
resolutions.  We follow the method used for cubic plane curves in
\cite{HartshorneAG} Proposition IV.4.21.

Let $\mathcal O := \mathcal O_{\mathbb P^5}$.  Then $\mathcal O$ has a
Frobenius endomorphism given on local sections by $x \mapsto x^p$.  If
the ideal sheaf $\mathcal I_X$ is locally generated by sections $f_1,
\dots, f_r$, then $F(\mathcal I_X)\mathcal O$ is locally generated by
$f_1^p, \dots f_r^p$.  Let $F(X)$ denote the subscheme of $\mathbb
P^5$ defined by the ideal sheaf $F(\mathcal I_X)\mathcal O$.  Then
$\mathcal I_{F(X)}\subset \mathcal I_X$, and the action of $F$ on
$\mathcal O_X$ factors as
\begin{equation}
\label{Frob.factor}
\mathcal O_X \rTo ^{F} \mathcal O_{F(X)}  \rOnto^{\pi} \mathcal O_X
\end{equation}
where $\pi$ is the canonical surjection.

Since $\mathbb P^5$ is a regular scheme, $\mathcal O$ is a flat
$F(\mathcal O)$-algebra (Kunz \cite{Kunz}).  So if one applies $F$ to
a locally free resolution of $\mathcal O_X$ one obtains a locally free
resolution of $\mathcal O_{F(X)}$.  Now the nonclassical (and
sometimes degenerate) Enriques surfaces $X \subset \mathbb P^5$ of
\S\S \ref{sect.rp2},\ref{sect.Enriques} have minimal free
resolutions similar to that of \eqref{free.res.P7}.  The factorization
\eqref{Frob.factor} lifts to the locally free resolutions and gives a
commutative diagram with exact rows:
\begin{footnotesize}
\begin{diagram}[w=2em]
0 & \rTo & \mathcal O(-6) & \rTo^d & \mathcal O(-5)^6 \oplus \mathcal
O(-6) & \rTo & \mathcal O(-4)^{15} & \rTo & \mathcal O(-3)^{10} &
\rTo^d & \mathcal O & \rOnto & \mathcal O_X \\
&& \dTo >F && \dTo >F && \dTo >F && \dTo >F && \dTo >F && \dTo
>F \\ 
0 & \rTo & \mathcal O(-12) & \rTo^{F(d)} & \mathcal O(-10)^6 \oplus
\mathcal O(-12) & \rTo & \mathcal O(-8)^{15} & \rTo &
\mathcal O(-6)^{10} & \rTo^{F(d)} & \mathcal O & \rOnto & \mathcal
O_{F(X)} \\
&& \dTo >g  && \dTo  && \dTo  && \dTo  && \dSame  && \dTo >\pi \\
0 & \rTo & \mathcal O(-6) & \rTo^d & \mathcal O(-5)^6 \oplus \mathcal
O(-6) & \rTo & \mathcal O(-4)^{15} & \rTo & \mathcal O(-3)^{10} &
\rTo^d & \mathcal O & \rOnto & \mathcal O_X
\end{diagram}%
\end{footnotesize}%
where the upper set of vertical maps are given by the action of
Frobenius, and the lower set of vertical maps are $\mathcal O$-linear
maps lifting $\pi$.  The action of Frobenius on the cohomology of
$\mathcal O_X$ can be computed using the hypercohomology of the
resolutions in the above diagram.  Only the terms farthest to the left
contribute to $H^1(\mathcal O_X)$.  So the action of Frobenius on
$H^1(\mathcal O_X)$ is the composite map
\[
H^5(\mathbb P^5, \mathcal O(-6)) \xrightarrow{F} H^5(\mathbb P^5,
\mathcal O(-12)) \xrightarrow{\,g\,} H^5(\mathbb P^5, \mathcal O(-6)).
\]
The action of Frobenius on the generator of the \v Cech $\check
H^5(\mathbb P^5, \mathcal O(-6))$, computed with respect to the
standard cover, is given by
\[
\frac 1{\prod x_i} \rMapsto \frac 1{(\prod x_i)^2} \rMapsto \coeff
\left( \frac g{(\prod x_i)^2}, \frac 1{\prod x_i} \right) =
\coeff(g,\prod x_i) \cdot \frac 1{\prod x_i}
\]
We conclude that the Hasse invariant of $X$ is the coefficient of
$\prod x_i$ in $g$.  If this coefficient is nonzero, then $X$ is a
$\mmu$-surface.  Otherwise it is an $\aalpha$-surface.

In the case of Reisner's monomial example $X(\Delta)$ in
\S\ref{sect.rp2} the minimal free resolution is $\mathbb
Z^6$-multigraded, and the minimal generators of each syzygy module
have different squarefree weights.  In such a situation, one can lift
the projection $\pi : R/F(I)R \to R/I$ to a chain map between minimal
free resolutions in which all the matrices are diagonal with nonzero
squarefree monomial entries.  The map $g : R(-12) \to R(-6)$ must
therefore be $\prod x_i$.  In particular the Hasse invariant is one.
This can also be verified with direct computations using {\em
Macaulay/Macaulay2} \cite{Macaulay} \cite{Macaulay2}.

{\em Macaulay/Macaulay2} can also be used in connection with the
methods of this chapter to produce examples of smooth Fano models of
unnodal nonclassical Enriques surface of both types $\mmu$ and
$\aalpha$ defined over the field with two elements.

The previous diagram may be globalized by taking the global resolution
\eqref{universal}, replacing $p^*\Omega_{\mathbb P^5}^3(0,-1)$ by its
Koszul free resolution, and applying Frobenius.  One obtains a diagram
of sheaves on $\mathbb P^5 \times \og_{20}$.
\begin{footnotesize}
\begin{diagram}
0 & \rTo & \mathcal O(-6,-1) & \rTo & \mathcal O(-5,-1)^6 \oplus
\mathcal O(-6,-2) & \rTo & \dotsb & \rTo & \mathcal O & \rOnto &
\mathcal O_Z \\ 
&& \dTo >F && \dTo >F &&  && \dTo >F && \dTo
>F \\ 
0 & \rTo & \mathcal O(-12,-2) & \rTo & \mathcal O(-10,-2)^6 \oplus
\mathcal O(-12,-4) & \rTo & \dotsb & \rTo & \mathcal O & \rOnto &
\mathcal O_{F(Z)} \\
&& \dTo >g  && \dTo  &&  && \dSame  && \dTo >\pi \\
0 & \rTo & \mathcal O(-6,-1) & \rTo & \mathcal O(-5,-1)^6 \oplus \mathcal
O(-6,-2) & \rTo & \dotsb & \rTo & \mathcal O & \rOnto & \mathcal O_Z
\end{diagram}%
\end{footnotesize}%
The coefficient of $\prod x_i$ in $g$ is now a section of
$\og_{20}(1)$.  We conclude:

\begin{proposition}
\label{Frob.locus}
The closure in $\og_{20}$ of the locus parametrizing Fano models of
unnodal $\aalpha$-surfaces is an $\SL_6$-invariant hyperplane section
of $\og_{20}$.
\end{proposition}

\

\end{document}